\documentclass[12pt,abstract]{scrartcl}
\usepackage[dvipdfmx]{graphicx}
\usepackage{amsmath,amsthm}
\usepackage[round]{natbib}
\newtheorem{proposition}{Proposition}
\usepackage{newtxtext,newtxmath}

\begin{document}
\title{A differential game analysis of R\&D in oligopoly with differentiated goods under general demand and cost functions: Bertrand vs. Cournot}

\author{Masahiko Hattori\thanks{%
mhattori@mail.doshisha.ac.jp} \\
Faculty of Economics, Hokkai-Gakuen University,\\
[.02cm] Toyohira-ku, Sapporo, Hokkaido, 062-8605, Japan,\\
[.01cm] \textrm{and} \\
[.1cm] Yasuhito Tanaka\thanks{%
yatanaka@mail.doshisha.ac.jp}\\
[.01cm] Faculty of Economics, Doshisha University,\\
Kamigyo-ku, Kyoto, 602-8580, Japan.\\
}

\date{}
\maketitle

\begin{abstract}
We study a dynamic oligopoly with differentiated goods by differential game approach under general demand and cost functions. We show that the steady state value of the R\&D investment by each firm is decreasing with respect to the number of firms, and the steady state value of the industry R\&D investment is increasing with respect to the number of firms. Also we show that if there is no spillover, whether the R\&D investment of each firm given the cost level in the memoryless closed-loop case is larger or smaller than that in the open-loop case depends on whether the strategic variables are strategic substitutes or strategic complements. Further we show that the memoryless closed-loop solution and the feedback solution (by the Hamilton-Jacobi-Bellman equation) are equivalent.
\end{abstract}

\begin{description}
\item[Keywords:] differential game; general demand and cost functions; R\&D in oligopoly with differentiated goods; open-loop; closed-loop; feedback
\end{description}

\begin{description}
\item[JEL Classification No.:] C73, D43, L13.
\end{description}

\section{Introduction}

In this paper we present an analysis of R\&D investment in a dynamic oligopoly model with differentiated goods by differential game approach. There are many studies of dynamic oligopoly by differential game theory, for example, \cite{cl0},  \cite{cl1}, \cite{cl3}, \cite{cl2}, \cite{cl4}, \cite{fu}, \cite{fu1}, and \cite{lam18}. Among them \cite{cl4} analyzed the problem of R\&D investment to cost-reducing activities in a Cournot oligopoly and a Bertrand oligopoly with differentiated goods. However, most of these studies including \cite{cl4} used a model of linear demand functions and quadratic or linear cost functions. These assumptions are very limited. We study the problem addressed by them in an oligopoly with general demand and cost functions.

In the next section we present a model and assumptions. In Section 3 we consider the steady state level of R\&D investment which is common to Bertrand and Cournot cases in the open-loop solution, the memoryless closed-loop and the feedback solution. In Section 4 we consider the open-loop solution of the R\&D investment in a Bertrand oligopoly. In Section 5 we consider the open-loop solution of the R\&D investment in a Cournot oligopoly, and compare the results of two cases in Section 6. We show the following results.
\begin{enumerate}
	\item The steady state value of the R\&D investment by each firm is decreasing with respect to the number of firms.
	\item The steady state value of the industry R\&D investment is increasing with respect to the number of firms.
\item The R\&D investment of each firm \emph{given the cost level} in the open-loop Bertrand oligopoly is larger than that in the open-loop Cournot oligopoly. 
\end{enumerate}

\cite{cl4} (also \cite{cl3} for an oligopoly with a homogeneous good) claim that the open-loop solution and the memoryless closed-loop solution coincide. However, as \cite{smr} point out, this claim is incorrect\footnote{Strictly speaking, \cite{smr} is a comment on \cite{cl5}. But \cite{cl3} and \cite{cl5} use the same model, and the model of \cite{cl4} is similar to it.} We present brief discussion about the memoryless closed-loop case in Section 7. We show the following results
\begin{enumerate}
\item  Suppose that there is no spillover effect of R\&D investment in a Bertrand oligopoly. If the strategic variables (prices) of the firms are strategic substitutes (or strategic complements), the R\&D investment of each firm given the cost level in the memoryless closed-loop case is larger (smaller) than that in the open-loop case.
\item  Suppose that there is no spillover effect of R\&D investment in a Cournot oligopoly. If the strategic variables (outputs) of the firms are strategic substitutes (or strategic complements), the R\&D investment of each firm given the cost level in the memoryless closed-loop case is larger (smaller) than that in the open-loop case.
\end{enumerate}

In Section 8 we examine the feedback solutions using the Hamilton-Jacobi-Bellman equation, and show that if there is no spillover effect of R\&D investment, the memoryless closed-loop solution and the feedback solution are equivalent both in the Bertrand oligopoly and the Cournot oligopoly.

\subsection*{Strategic substitutability and strategic complementarity}

We assume that the goods of the firms are \emph{substitutes}, not complements. This means that a rise in the price of one good increases the demands for other goods, and an increase in the output of one good lowers the prices of other goods. However, the strategic variables (outputs or prices) of the firms may be strategic substitutes or strategic complements. 

In the Cournot oligopoly, if the reaction of a firm's output to an increase in the output of another firm is negative (or positive), the outputs of firms are strategic substitutes (or strategic complements).

Note that if inverse demand (and direct) functions are linear, the outputs are strategic substitutes. 

On the other hand, in the Bertrand oligopoly, if the reaction of the price of a firm's good to a rise in the price of another firm' good is negative (or positive), the prices of the goods of firms are strategic substitutes (or strategic complements).

 Note that if direct demand (and inverse demand) functions are linear, the prices are \emph{strategic complements}. 

\section{The model}

Consider an oligopoly with $n$ firms in which at any $t\in [0, \infty)$  they produce differentiated goods to maximize their discounted profits. The goods are substitutes. The firms are called Firms 1, 2, $\dots$, $n$. Let $q_i(t)$ be the output of Firm $i,\ i\in \{1, 2, \dots, n\}$, $p_i(t)$ be the price of the good of Firm $i,\ i\in \{1, 2, \dots, n\}$, at $t$. 
The utility of a representative consumer is
\[u(q_1(t),q_2(t),\dots,q_n(t))+x(t).\]
$x(t)$ is the consumption of a numeraire good. Let $y(t)$ be his income. Then, the utility maximization problem is
\[\max_{q_1(t),q_2(t),\dots,q_n(t)}[u(q_1(t),q_2(t),\dots,q_n(t))+x(t)]\]
subject to
\[\sum_{i=1}^np_i(t)q_i(t)+x(t)=y(t).\]
The conditions for utility maximization are
\[\frac{\partial u(q_1(t),q_2(t),\dots,q_n(t))}{\partial q_i(t)}=p_i(t),\ i\in \{1, 2, \dots, n\}.\]
From them the inverse demand functions are obtained as follows.
\[p_i(t)=p_i(q_1(t),q_2(t),\dots,q_n(t)),\ i\in \{1, 2, \dots, n\}.\]
We have $\frac{\partial p_i(q_1(t),q_2(t),\dots,q_n(t))}{\partial q_i(t)}<0,\ i\in \{1, 2, \dots, n\}$, and since the goods are substitutes
\[\frac{\partial p_i(q_1(t),q_2(t),\dots,q_n(t))}{\partial q_j(t)}<0.\]
If the outputs of the firms are strategic substitutes,
\[\frac{\partial p_i(q_1(t),q_2(t),\dots,q_n(t))}{\partial q_j(t)}+\frac{\partial^2 p_i(q_1(t),q_2(t),\dots,q_n(t))}{\partial q_i(t)\partial q_j(t)}q_i(t)<0,\ j\neq i.\]
If they are strategic complements,
\[\frac{\partial p_i(q_1(t),q_2(t),\dots,q_n(t))}{\partial q_j(t)}+\frac{\partial^2 p_i(q_1(t),q_2(t),\dots,q_n(t))}{\partial q_i(t)\partial q_j(t)}q_i(t)>0,\ j\neq i.\]

The production cost of Firm $i$ is
\[C(m_i(t),q_i(t)),\ i\in \{1, 2,\dots,n\}.\]
$m_i(t)$ is a parameter which represents the current cost of Firm $i$. Denote $C(m_i(t),q_i(t))$ by $C_i$.  It satisfies
\[\frac{\partial C_i}{\partial q_i(t)}>0,\ \frac{\partial C_i}{\partial m_i(t)}>0,\ \frac{\partial^2 C_i}{\partial q_i(t)\partial m_i(t)}>0.\]

Inverting the inverse demand functions, the direct demand functions are obtained as follows.
\[q_i(t)=q(p_1(t),p_2(t),\dots,p_n(t)),\ i\in \{1, 2, \dots, n\}.\]
We have $\frac{\partial q_i(p_1(t),p_2(t),\dots,p_n(t))}{\partial p_i(t)}<0,\ i\in \{1, 2, \dots, n\}$, and since the goods are substitutes
\[\frac{\partial q_i(p_1(t),p_2(t),\dots,p_n(t))}{\partial p_j(t)}>0.\]
If the prices of the goods of the firms are strategic substitutes,
\[\frac{\partial q_i}{\partial p_j(t)}+\left(p_i(t)-\frac{\partial C_i}{\partial q_i(t)}\right)\frac{\partial^2 q_i}{\partial p_i(t)\partial p_j(t)}<0,\ j\neq i,\]
If they are strategic complements,
\[\frac{\partial q_i}{\partial p_j(t)}+\left(p_i(t)-\frac{\partial C_i}{\partial q_i(t)}\right)\frac{\partial^2 q_i}{\partial p_i(t)\partial p_j(t)}>0,\ j\neq i.\]

Let $k_i(t)$ be the R\&D investment by Firm $i$. The moving of $m_i(t)$ is governed by 
\begin{equation}
\frac{dm_i(t)}{dt}=m_i(t)[-\Gamma(k_i(t),K_{-i}(t))+\delta],\label{m1}
\end{equation}
where
\[\Gamma(k_i(t),K_{-i}(t))>0,\ K_{-i}=\sum_{j\neq i}k_j(t).\]
$\delta\in (0,1)$ is a constant depreciation rate. Denote $\Gamma(k_i(t),K_{-i}(t))$ by $\Gamma_i$. We assume that $\Gamma_i$ is strictly increasing and concave, that is,
\[\frac{\partial \Gamma_i}{\partial k_i(t)}>0,\ \frac{\partial \Gamma_i}{\partial K_{-i}(t)}\geq 0,\]
and
\[\frac{\partial^2 \Gamma_i}{\partial k_i(t)^2}\leq 0,\ \frac{\partial^2 \Gamma_i}{\partial K_{-i}(t)^2}\leq 0,\ \frac{\partial^2 \Gamma_i}{\partial k_i(t)\partial K_{-i}(t)}\leq 0.\]
Also we assume
\[\left|\frac{\partial \Gamma(k_i(t),K_{-i}(t))}{\partial k_i(t)}\right|>\left|\frac{\partial \Gamma(k_i(t),K_{-i}(t))}{\partial K_{-i}(t)}\right|.\]
This means that the direct effect of R\&D investment is larger than the spillover effect.

The R\&D cost of Firm $i$ is
\[\gamma(k_i(t)),\ i\in \{1, 2,\dots,n\}.\]
We assume that it is strictly increasing and strictly convex, that is, $\gamma'(k_i(t))>0$ and $\gamma''(k_i(t))>0$.

\section{The steady state R\&D investment}

Let $k^*$ be the steady state value of $k_i(t)$. At the steady state the following equation holds.
\begin{equation}
\Gamma(k^*,(n-1)k^*)=-\delta.\label{st}
\end{equation}
From this we obtain
\begin{equation}
\frac{dk^*}{dn}=-\frac{\frac{\partial \Gamma(k^*,(n-1)k^*)}{\partial K_{-i}(t)}k^*}{\frac{\partial \Gamma(k^*,(n-1)k^*)}{\partial k_i(t)}+(n-1)\frac{\partial \Gamma(k^*,(n-1)k^*)}{\partial K_{-i}(t)}}\leq 0.\label{k1}
\end{equation}
Also
\begin{equation}
\frac{dnk^*}{dn}=k^*+n\frac{dk^*}{dn}=\frac{\left(\frac{\partial \Gamma(k^*,(n-1)k^*)}{\partial k_{i}(t)}-\frac{\partial \Gamma(k^*,(n-1)k^*)}{\partial K_{-i}(t)}\right)k^*}{\frac{\partial \Gamma(k^*,(n-1)k^*)}{\partial k_i(t)}+(n-1)\frac{\partial \Gamma(k^*,(n-1)k^*)}{\partial K_{-i}(t)}}>0.\label{k2}
\end{equation}
They are because $\frac{\partial \Gamma(k^*,(n-1)k^*)}{\partial k_{i}(t)}>0$, $\frac{\partial \Gamma(k^*,(n-1)k^*)}{\partial K_{-i}(t)}\geq 0$, $\left|\frac{\partial \Gamma(k^*,(n-1)k^*)}{\partial k_{i}(t)}\right|>\left|\frac{\partial \Gamma(k^*,(n-1)k^*)}{\partial K_{-i}(t)}\right|$.

These results mean that the steady state value of the R\&D investment by each firm is decreasing with respect to the number of firms, while the total R\&D investment is increasing with respect to the number of firms.

Summarizing the results in the following proposition.
\begin{proposition}
\begin{enumerate}
	\item The steady state value of the R\&D investment by each firm is decreasing with respect to the number of firms.
	\item The steady state value of the industry R\&D investment is increasing with respect to the number of firms.
\end{enumerate}
\end{proposition}

Note that these conclusions hold in both the Bertrand oligopoly and the Cournot oligopoly and in the open-loop case, the memoryless closed-loop case and the feedback case because (\ref{st}) holds in all cases.

\section{R\&D in a dynamic oligopoly: Bertrand competition}

We seek to the solution of the open-loop approach in the Bertrand oligopoly.  The instantaneous profit of Firm $i$ is written as
\[q_i(p_1(t),p_2(t),\dots,p_n(t))p_i(t)-C(m_i(t),q_i(t))-\gamma(k_i(t)).\]
The objective of Firm $i$ is 
\[\max_{p_i(t), k_i(t)}\int_{0}^{\infty}e^{-\rho t}[q_i(p_1(t),p_2(t),\dots,p_n(t))p_i(t)-C(m_i(t),q_i(t))-\gamma(k_i(t))]dt,\]
subject to (\ref{m1}).

The present value Hamiltonian function for Firm $i,\ i\in \{1, 2, \dots, n\}$, is
\begin{align*}
\mathcal{H}_i(p_i(t), k_i(t))=&e^{-\rho t}\big\{q_i(p_1(t),p_2(t),\dots,p_n(t))p_i(t)-C(m_i(t),q_i(t))-\gamma(k_i(t))\\
&+\lambda_{ii}(t)m_i(t)[-\Gamma(k_i(t),K_{-i}(t))+\delta]+\sum_{j\neq i}\lambda_{ij}(t)m_j(t)[-\Gamma(k_j(t),K_{-j}(t))+\delta]\big\},
\end{align*}
where
\[K_{-j}(t)=\sum_{l\neq j}k_l(t).\]
This includes $k_i(t)$.

The current value Hamiltonian function for Firm $i,\ i\in \{1, 2, \dots, n\}$, is
\begin{align*}
\hat{\mathcal{H}}_i(p_i(t), k_i(t))=&e^{\rho t}\mathcal{H}_i=q_i(p_1(t),p_2(t),\dots,p_n(t))p_i(t)-C(m_i(t),q_i(t))-\gamma(k_i(t))\\
&+\lambda_{ii}(t)m_i(t)[-\Gamma(k_i(t),K_{-i}(t))+\delta]+\sum_{j\neq i}\lambda_{ij}(t)m_j(t)[-\Gamma(k_j(t),K_{-j}(t))+\delta].
\end{align*}
Let
\[\mu_i(t)=e^{-\rho t}\lambda_i.\]
$\mu_i(t)$ is the costate variable. Denote $\hat{\mathcal{H}}_i(p_i(t), k_i(t))$ by $\hat{\mathcal{H}}_i$. 

The first order conditions for Firm $i$ are
\begin{align}
\frac{\partial \hat{\mathcal{H}}_i}{\partial p_i(t)}=&q_i(p_1(t),p_2(t),\dots,p_n(t)) \label{q}\\
&+\left(p_i(t)-\frac{\partial C(m_i(t),q_i(t))}{\partial q_i(t)}\right)\frac{\partial q_i(p_1(t),p_2(t),\dots,p_n(t)}{\partial p_i(t)}=0,\notag
\end{align}
and
\begin{equation}
\frac{\partial \hat{\mathcal{H}}_i}{\partial k_i(t)}=-\gamma'(k_i(t))-\lambda_{ii}(t)\frac{\partial \Gamma_i}{\partial k_i(t)}m_i(t)-\sum_{j\neq i}\lambda_{ij}(t)\frac{\partial \Gamma_i}{\partial K_{-i}(t)}m_j(t)=0.\label{k}
\end{equation}
The second order condition for Firm $i$ about the price choice is
\begin{align*}
\frac{\partial \hat{\mathcal{H}}^2_i}{\partial p_i(t)^2}=&2\frac{\partial q_i(p_1(t),p_2(t),\dots,p_n(t)}{\partial p_i(t)}+\left(p_i(t)-\frac{\partial C(m_i(t),q_i(t))}{\partial q_i(t)}\right)\frac{\partial^2 q_i(p_1(t),p_2(t),\dots,p_n(t)}{\partial p_i(t)^2}\\
&-\frac{\partial^2 C(m_i(t),q_i(t))}{\partial q_i(t)^2}\left(\frac{\partial q_i(p_1(t),p_2(t),\dots,p_n(t)}{\partial p_i(t)}\right)^2<0.
\end{align*}
Its second order condition about the R\&D investment choice is
\begin{equation}
\frac{\partial^2 \hat{\mathcal{H}}_i}{\partial k_i(t)^2}=-\gamma''(k_i(t))-\lambda_{ii}(t)\frac{\partial^2 \Gamma_i}{\partial k_i(t)^2}m_i(t)-\sum_{j\neq i}\lambda_{ij}(t)\frac{\partial^2 \Gamma_j}{\partial K_{-i}(t)^2}m_j(t)<0.\label{sc1}
\end{equation}
The adjoint conditions are
\begin{equation}
-\frac{\partial \hat{\mathcal{H}}_i}{\partial m_i(t)}=\frac{\partial \lambda_{ii}(t)}{\partial t}-\rho \lambda_{ii}(t),\ i\in \{1, 2, \dots,n\},\label{a1}
\end{equation}
and
\begin{equation}
-\frac{\partial \hat{\mathcal{H}}_i}{\partial m_j(t)}=\frac{\partial \lambda_{ij}(t)}{\partial t}-\rho \lambda_{ij}(t),\ j\neq i.\label{a2}
\end{equation}
We have
\begin{equation}
\frac{\partial \hat{\mathcal{H}}_i}{\partial m_i(t)}=-\frac{\partial C_i}{\partial m_i(t)}+\lambda_{ii}(t)[-\Gamma(k_i(t),K_{-i}(t))+\delta],\label{b1}
\end{equation}
\begin{equation}
\frac{\partial \hat{\mathcal{H}}_i}{\partial m_j(t)}=\lambda_{ij}(t)[-\Gamma(k_j(t),K_{-j}(t)+\delta],\label{b2}
\end{equation}
At the steady state
\[\frac{dm_i(t)}{dt}=m_i(t)[-\Gamma(k_i(t),K_{-i}(t))+\delta]=0,\]
$\frac{\partial \lambda_{ii}}{\partial t}=0$ and $\frac{\partial \lambda_{ij}}{\partial t}=0,\ i\in \{1, 2, \dots, n\},\ j\neq i$. By symmetry of the oligopoly we can assume $\lambda_{ii}(t)=\lambda_{jj}(t)$ for $j\neq i$, $\lambda_{ij}(t)=\lambda_{il}(t)=\lambda_{ji}(t)$ for $j, l\neq i$, $\frac{\partial k_j(t)}{\partial m_i(t)}=\frac{\partial k_l(t)}{\partial m_j(t)}$ for $l\neq j$, $m_i(t)=m_j(t)$ for $j\neq i$, and so on. Denote the steady state values of $\lambda_{ii}$, $\lambda_{ij}$, $p_i(t)$, $k_i(t)$ and $m_i(t)$ by $\lambda_{own}$, $\lambda_{other}$, $p^{*}_B$, $k^{*}_B$ and $m^{*}_B$. Then, (\ref{b1}) and (\ref{b2}) are reduced to
\begin{equation*}
\frac{\partial \hat{\mathcal{H}}_i}{\partial m_i(t)}=-\frac{\partial C_i}{\partial m_i(t)},
\end{equation*}
\begin{equation*}
\frac{\partial \hat{\mathcal{H}}_i}{\partial m_j(t)}=0,
\end{equation*}
and (\ref{a1}) and (\ref{a2}) are rewritten as
\[\frac{\partial C_i}{\partial m_i(t)}=-\rho \lambda_{own},\]
and
\[\lambda_{other}=0.\]
Since $\frac{\partial C_i}{\partial m_i(t)}<0$, we have $\lambda_{own}>0$. The first order condition for the choice of $k_i(t)$, (\ref{k}), is reduced to
\[-\gamma'(k_i(t))-\lambda_{own}\frac{\partial \Gamma_i}{\partial k_i(t)}m^{*}_B=0.\]
This means
\begin{equation}
\frac{\partial C(m^{*}_B,q^{*}_B)}{\partial m_i(t)}\frac{\partial \Gamma(k^{*}_B,(n-1)k^{*}_B)}{\partial k_i(t)}=\frac{\rho}{m^{*}_B}\gamma'(k^{*}_B).\label{x1}
\end{equation}

\subsection*{\bfseries Linear and quadratic example}

According to \cite{cl4}, assume that the direct demand functions are
\[q_i(t)=\frac{a}{\varphi}-\frac{(\varphi-s)p_i(t)}{(1-s)\varphi}+\frac{s}{(1-s)\varphi}\sum_{j\neq i}p_j(t),\]
where $\varphi=1+(n-1)s.$
The production  cost of Firm $i,\ i\in \{1, 2, \dots, n\}$, is
\[C(m_i(t),q_i(t))=m_i(t)q_i(t),\]
the R\&D cost of Firm $i,\ i\in \{1, 2, \dots, n\}$, is
\[\gamma(k_i(t))=b[k_i(t)]^2,\ b>0.\]
The moving of $m_i(t)$ is governed by
\[\frac{dm_i(t)}{dt}=m_i(t)[-\Gamma(k_i(t),K_{-i}(t))+\delta]=m_i(t)[-k_i(t)-\beta K_{-i}(t)+\delta],\ \beta<1.\](\ref{k1}) and (\ref{k2}) are reduced to
\[\frac{dk^{*}_B}{dn}=-\frac{\beta k^{*}_B}{1+(n-1)\beta}<0,\]
and
\[\frac{dnk^{*}_B}{dn}=k^{*}_B+n\frac{dk^{*}_B}{dn}=\frac{(1-\beta)k^{*}_B}{1+(n-1)\beta}>0.\]
From (\ref{x1})
\[q^{*}_B=\frac{2\rho}{m^{*}_B}bk^{*}_B.\]
Since
\[p^{*}_B=\frac{a(1-s)+m^{*}_B[1+(n-2)s]}{2+(n-3)s},\ q^{*}_B=\frac{(a-m^{*}_B)[1+(n-2)s]}{[2+(n-3)s][1+(n-1)s]}\]
in this example, we have
\[k^{*}_B\big|_{\mathrm{given}\ m^{*}_B}=\frac{m^{*}_B(a-m^{*}_B)[1+(n-2)s]}{2b\rho[2+(n-3)s][1+(n-1)s]}.\]

\section{R\&D in a dynamic oligopoly: Cournot competition}

We seek to the solution of the open-loop approach in the Cournot oligopoly.  The instantaneous profit of Firm $i$ is written as
\[p_i(q_1(t),q_2(t),\dots,q_n(t))q_i(t)-C(m_i(t),q_i(t))-\gamma(k_i(t)).\]
The objective of Firm $i$ is 
\[\max_{q_i(t), k_i(t)}\int_{0}^{\infty}e^{-\rho t}[p_i(q_1(t),q_2(t),\dots,q_n(t))q_i(t)-C(m_i(t),q_i(t))-\gamma(k_i(t))]dt,\]
subject to (\ref{m1}).

The present value Hamiltonian function for Firm $i,\ i\in \{1, 2, \dots, n\}$ is
\begin{align*}
\mathcal{H}_i(p_i(t), k_i(t))=&e^{-\rho t}\big\{p_i(q_1(t),q_2(t),\dots,q_n(t))q_i(t)-C(m_i(t),q_i(t))-\gamma(k_i(t))\\
&+\lambda_{ii}(t)m_i(t)[-\Gamma(k_i(t),K_{-i}(t))+\delta]+\sum_{j\neq i}\lambda_{ij}(t)m_j(t)[-\Gamma(k_j(t),K_{-j}(t))+\delta]\big\},
\end{align*}
where
\[K_{-j}(t)=\sum_{l\neq j}k_l(t).\]

The current value Hamiltonian function for Firm $i,\ i\in \{1, 2, \dots, n\}$ is
\begin{align*}
\hat{\mathcal{H}}_i(p_i(t), k_i(t))=&e^{\rho t}\mathcal{H}_i=p_i(q_1(t),q_2(t),\dots,q_n(t))q_i(t)-C(m_i(t),q_i(t))-\gamma(k_i(t))\\
&+\lambda_{ii}(t)m_i(t)[-\Gamma(k_i(t),K_{-i}(t))+\delta]+\sum_{j\neq i}\lambda_{ij}(t)m_j(t)[-\Gamma(k_j(t),K_{-j}(t))+\delta].
\end{align*}
Let
\[\mu_i(t)=e^{-\rho t}\lambda_i.\]
$\mu_i(t)$ is the costate variable. Denote $\hat{\mathcal{H}}_i(q_i(t), k_i(t))$ by $\hat{\mathcal{H}}_i$. 

The first order conditions for Firm $i$ are
\begin{align}
\frac{\partial \hat{\mathcal{H}}_i}{\partial q_i(t)}=&p_i(q_1(t),q_2(t),\dots,q_n(t))+\frac{\partial p_i(q_1(t),q_2(t),\dots,q_n(t))}{\partial q_i(t)}q_i(t)\label{cq}\\
&-\frac{\partial C(m_i(t),q_i(t))}{\partial q_i(t)}=0,\notag
\end{align}
and
\begin{equation*}
\frac{\partial \hat{\mathcal{H}}_i}{\partial k_i(t)}=-\gamma'(k_i(t))-\lambda_{ii}(t)\frac{\partial \Gamma_i}{\partial k_i(t)}m_i(t)-\sum_{j\neq i}\lambda_{ij}(t)\frac{\partial \Gamma_i}{\partial K_{-i}(t)}m_j(t)=0.
\end{equation*}
The second order condition for Firm $i$ about the output choice is
\[\frac{\partial \hat{\mathcal{H}}^2_i}{\partial q_i(t)^2}=2\frac{\partial p_i(q_1(t),q_2(t),\dots,q_n(t))}{\partial q_i(t)}+\frac{\partial^2 p_i(q_1(t),q_2(t),\dots,q_n(t))}{\partial q_i(t)^2}q_i(t)-\frac{\partial^2 C(m_i(t),q_i(t))^2}{\partial q_i(t)}<0.\]
Its second order condition about the R\&D investment choice is
\begin{equation}
\frac{\partial^2 \hat{\mathcal{H}}_i}{\partial k_i(t)^2}=-\gamma''(k_i(t))-\lambda_{ii}(t)\frac{\partial^2 \Gamma_i}{\partial k_i(t)^2}m_i(t)-\sum_{j\neq i}\lambda_{ij}(t)\frac{\partial^2 \Gamma_j}{\partial K_{-i}(t)^2}m_j(t)<0.\label{sc2}
\end{equation}
The adjoint conditions are the same those in the Bertrand oligopoly. Denote the steady state values of $\lambda_{ii}$, $\lambda_{ij}$, $q_i(t)$, $k_i(t)$ and $m_i(t)$ by $\lambda_{own}$, $\lambda_{other}$, $q^{*}_C$, $k^{*}_C$ and $m^{*}_C$. Similarly to the Bertrand oligopoly we obtain
\begin{equation}
\frac{\partial C(m^{*}_C,q^{*}_C)}{\partial m_i(t)}\frac{\partial \Gamma(k^{*}_C,(n-1)k^{*}_C)}{\partial k_i(t)}=\frac{\rho}{m^{*}_C}\gamma'(k^{*}_C).\label{cx1}
\end{equation}

\subsection*{\bfseries Linear and quadratic example}

Assume that the inverse demand functions are
\[p_i(t)=a-q_i(t)-\sum_{j\neq i}^nsq_j(t),\ 0<s<1.\]
They are derived from the direct demand functions in the previous example.

The production  cost of Firm $i,\ i\in \{1, 2, \dots, n\}$ is
\[C(m_i(t),q_i(t))=m_i(t)q_i(t),\]
the R\&D investment cost of Firm $i,\ i\in \{1, 2, \dots, n\}$, is
\[\gamma(k_i(t))=b[k_i(t)]^2,\ b>0.\]
The moving of $m_i(t)$ is governed by
\[\frac{dm_i(t)}{dt}=m_i(t)[-\Gamma(k_i(t),K_{-i}(t))+\delta]=m_i(t)[-k_i(t)-\beta K_{-i}(t)+\delta],\ \beta<1.\]From (\ref{cx1})
\[q^{*}_C=\frac{2\rho}{m^{*}_C}bk^{*}_C.\]
Since $q^{*}_C=\frac{a-m^{*}_C}{2+(n-1)s}$ in this example,
\[k^{*}_C\big|_{\mathrm{given}\ m^{*}_C}=\frac{m^{*}_C(a-m^{*}_C)}{2b\rho[2+(n-1)s]}.\]

\section{Comparison of Bertrand and Cournot given the cost level}

In the linear and quadratic example, as \cite{cl4} shows, we have
\[k^{*}_B>k^{*}_C.\]
We examine the general case. From (\ref{x1}) and (\ref{cx1}) $k^{*}_B$ (or $k^{*}_C$, the same hereinafter) given the cost level $m^{*}_B=m^{*}_C$ is obtained by the following equation.
\begin{equation}
-\frac{\rho}{m^{*}_B}\gamma'(k^{*}_B)+\frac{\partial C(m^{*}_B,q^{*}_B)}{\partial m_i(t)}\frac{\partial \Gamma(k^{*}_B,(n-1)k^{*}_B)}{\partial k_i(t)}=0.\label{ap1}
\end{equation}
By the second order conditions about the R\&D investment, (\ref{sc1}) and (\ref{sc2}), the left-hand side of (\ref{ap1}) is decreasing with respect to $k^{*}_B$. On the other hand, since $\frac{\partial^2 C(m^{*}_B,q^{*}_B)}{\partial q_i(t) \partial m_i(t)}>0$, $\frac{\partial C(m^{*}_B,q^{*}_B)}{\partial m_i(t)}$ is increasing with respect to $q^{*}_B$. Therefore, the larger is $q^{*}_B$, the larger is $k^{*}_B$. 

The first order condition for the price choice in the Bertrand oligopoly is
\begin{equation}
q_i(p^{*}_B,p^{*}_B,\dots,p^{*}_B)+\left(p^{*}_B-\frac{\partial C(m_i(t),q_i(t))}{\partial q_i(t)}\right)\frac{\partial q_i(p^{*}_B,p^{*}_B,\dots,p^{*}_B)}{\partial p_i(t)}=0.\label{ap3}
\end{equation}
The first order condition for the output choice in the Cournot oligopoly is
\begin{equation}
p_i(q^{*}_C, q^{*}_C, \dots, q^{*}_C)+\frac{\partial p_i(q^{*}_C, q^{*}_C, \dots, q^{*}_C)}{\partial q_i(t)}q^{*}_C-\frac{\partial C(m_i(t),q_i(t))}{\partial q_i(t)}=0.\label{ap4}
\end{equation}
Assume $q_i(t)=q^{*}_B$ and $p_i(t)=p^{*}_B$ for all $i\in \{1, 2, \dots, n\}$. From (\ref{ap3})
\[q^{*}_B=-\left(p^{*}_B-\frac{\partial C(m_i(t),q_i(t))}{\partial q_i(t)}\right)\frac{\partial q_i(p^{*}_B,p^{*}_B,\dots,p^{*}_B)}{\partial p_i(t)}.\]
Substituting this into the left-hand side of (\ref{ap4}) (assuming $q^{*}_C=q^{*}_B$) yields
\begin{equation}
\left(p^{*}_B-\frac{\partial C(m_i(t),q_i(t))}{\partial q_i(t)}\right)\left(1-\frac{\partial q_i(p^{*}_B,p^{*}_B,\dots,p^{*}_B)}{\partial p_i(t)}\frac{\partial p_i(q^{*}_B, q^{*}_B, \dots, q^{*}_B)}{\partial q_i(t)}\right).\label{ap5}
\end{equation}
Since $p^{*}_B-\frac{\partial C(m_i(t),q_i(t))}{\partial q_i(t)}>0$, and from (\ref{1-a}) in Appendix 1
\begin{equation}
\frac{\partial q_i(p^{*}_B,p^{*}_B,\dots,p^{*}_B)}{\partial p_i(t)}\frac{\partial p_i(q^{*}_B, q^{*}_B, \dots, q^{*}_B)}{\partial q_i(t)}>1,\label{ee1}
\end{equation}
(\ref{ap5}) is negative, and then the output of each firm in the Bertrand oligopoly is larger than that in the Cournot oligopoly. Thus, we obtain the following proposition.
\begin{proposition}
The R\&D investment of each firm given the cost level in the open-loop Bertrand oligopoly is larger than that in the open-loop Cournot oligopoly. 
\end{proposition}

\section{Memoryless closed-loop solution without spillover}

\subsection{Bertrand oligopoly}

We seek to the solution of the memoryless closed-loop approach in the Bertrand oligopoly. For simplicity, we assume 
\[\frac{\partial \Gamma_i}{\partial K_{-i}(t)}=0,\]
that is, there is no spillover effect of R\&D investment. The first order conditions for Firm $i$ are
\begin{align}
\frac{\partial \hat{\mathcal{H}}_i}{\partial p_i(t)}=&q_i(p_1(t),p_2(t),\dots,p_n(t)) \label{cbq}\\
&+\left(p_i(t)-\frac{\partial C(m_i(t),q_i(t))}{\partial q_i(t)}\right)\frac{\partial q_i(p_1(t),p_2(t),\dots,p_n(t)}{\partial p_i(t)}=0,\notag
\end{align}
and
\begin{equation}
\frac{\partial \hat{\mathcal{H}}_i}{\partial k_i(t)}=-\gamma'(k_i(t))-\lambda_{ii}(t)\frac{\partial \Gamma_i}{\partial k_i(t)}m_i(t)=0.\label{cbk}
\end{equation}
The adjoint conditions are
\begin{align}
&-\frac{\partial \hat{\mathcal{H}}_i}{\partial m_i(t)}-\sum_{j\neq i}\frac{\partial \hat{\mathcal{H}}_i}{\partial k_j(t)}\frac{\partial k_j(t)}{\partial m_i(t)}-\sum_{j\neq i}\frac{\partial \hat{\mathcal{H}}_i}{\partial p_j(t)}\frac{\partial p_j(t)}{\partial m_i(t)}\label{a21}\\
=&\frac{\partial \lambda_{ii}(t)}{\partial t}-\rho \lambda_{ii}(t),\ i\in \{1, 2, \dots,n\}, \notag
\end{align}
and
\begin{align*}
&-\frac{\partial \hat{\mathcal{H}}_i}{\partial m_j(t)}-\sum_{l\neq i,j}\frac{\partial \hat{\mathcal{H}}_i}{\partial k_l(t)}\frac{\partial k_l(t)}{\partial m_j(t)}-\frac{\partial \hat{\mathcal{H}}_i}{\partial k_j(t)}\frac{\partial k_j(t)}{\partial m_j(t)}\\
&-\sum_{l\neq i,j}\frac{\partial \hat{\mathcal{H}}_i}{\partial p_l(t)}\frac{\partial p_l(t)}{\partial m_j(t)}-\frac{\partial \hat{\mathcal{H}}_i}{\partial p_j(t)}\frac{\partial p_j(t)}{\partial m_j(t)}=\frac{\partial \lambda_{ij}(t)}{\partial t}-\rho \lambda_{ij}(t),\ j\neq i.\notag
\end{align*}
We have
\begin{equation}
\frac{\partial \hat{\mathcal{H}}_i}{\partial m_i(t)}=-\frac{\partial C_i}{\partial m_i(t)}+\lambda_{ii}(t)[-\Gamma(k_i(t),K_{-i}(t))+\delta],\label{b11}
\end{equation}
\begin{equation}
\frac{\partial \hat{\mathcal{H}}_i}{\partial m_j(t)}=\lambda_{ij}(t)[-\Gamma(k_j(t),K_{-j}(t))+\delta], \label{b21}
\end{equation}
\begin{equation*}
\frac{\partial \hat{\mathcal{H}}_i}{\partial k_j(t)}=-\lambda_{ij}(t)\frac{\partial \Gamma_j}{\partial k_j(t)}m_j(t),
\end{equation*}
\begin{equation*}
\frac{\partial k_j(t)}{\partial m_i(t)}=0,
\end{equation*}
\begin{equation*}
\frac{\partial k_j(t)}{\partial m_j(t)}=-\frac{\lambda_{jj}(t)}{\gamma''(k_i(t))+\lambda_{jj}(t)\frac{\partial^2 \Gamma_j}{\partial k_i(t^2)}m_j(t)}\frac{\partial \Gamma_j}{\partial k_j(t)},
\end{equation*}
and
\[\frac{\partial \hat{\mathcal{H}}_i}{\partial p_j(t)}=\frac{\partial q_i}{\partial p_j(t)}p_i(t)\]
$\frac{\partial p_j(t)}{\partial m_i(t)}$ is obtained by (\ref{2ap3}) in Appendix 2. If the prices of goods of the firms are strategic substitutes $\left(\frac{\partial q_j}{\partial p_i(t)}+\left(p_j(t)-\frac{\partial C(m_j(t),q_j(t))}{\partial q_j(t)}\right)\frac{\partial^2 q_j}{\partial p_i(t)\partial p_j(t)}<0\right)$, $\frac{\partial p_j(t)}{\partial m_i(t)}<0$, and if they are strategic complements $\left(\frac{\partial q_j}{\partial q_i(t)}+\left(p_j(t)-\frac{\partial C(m_j(t),q_j(t))}{\partial q_j(t)}\right)\frac{\partial^2 q_j}{\partial p_i(t)\partial p_j(t)}>0\right)$, $\frac{\partial p_j(t)}{\partial m_i(t)}>0$.

At the steady state we have
\[\frac{dm_i(t)}{dt}=m_i(t)[-\Gamma(k_i(t), K_{-i}(t))+\delta]=0,\]
$\frac{\partial \lambda_{ii}}{\partial t}=0$ and $\frac{\partial \lambda_{ij}}{\partial t}=0,\ i\in \{1, 2, \dots, n\},\ j\neq i$. Then, (\ref{b11}) and (\ref{b21}) are reduced to
\begin{equation*}
\frac{\partial \hat{\mathcal{H}}_i}{\partial m_i(t)}=-\frac{\partial C_i}{\partial m_i(t)},
\end{equation*}
\begin{equation*}
\frac{\partial \hat{\mathcal{H}}_i}{\partial m_j(t)}=0.
\end{equation*}
By symmetry of the oligopoly we can assume $\lambda_{ii}(t)=\lambda_{jj}(t)$ for $j\neq i$, $\lambda_{ij}(t)=\lambda_{il}(t)=\lambda_{ji}(t)$ for $j, l\neq i$, $\frac{\partial k_j(t)}{\partial m_i(t)}=\frac{\partial k_l(t)}{\partial m_j(t)}$ for $l\neq j$, $m_i(t)=m_j(t)$ for $j\neq i$, $p_j(t)=p_i(t)$, and so on.  Denote the steady state values of $\lambda_{ii}$, $\lambda_{ij}$, $p_i(t)$ and $k_i(t)$ by $\lambda_{own}$, $\lambda_{other}$, $p^{**}_B$ and $k^{**}_B$. Then, (\ref{a21}) is rewritten as
\[\frac{\partial C_i}{\partial m_i(t)}-(n-1)\frac{\partial q_i}{\partial p_j(t)}p^{**}_B\frac{\partial p_j(t)}{\partial m_i(t)}=-\rho \lambda_{own}.\]
The first order condition for the choice of $k_i(t)$, (\ref{cbk}), is reduced to
\[-\gamma'(k^{**}_B)-\lambda_{own}\frac{\partial \Gamma_i}{\partial k_i(t)}m_i(t)=0.\]
This means
\begin{align}
&\left[\frac{\partial C(m_i(t),q_i(t))}{\partial m_i(t)}-(n-1)\frac{\partial q_i}{\partial p_j(t)}p^{**}_B\frac{\partial p_j(t)}{\partial m_i(t)}\right]\frac{\partial \Gamma(k^{**}_B,(n-1)k^{**}_B)}{\partial k_i(t)}\label{x2}\\
&=\frac{\rho}{m_i(t)}\gamma'(k^{**}_B).\notag
\end{align}
The first order condition for the price choice in the memoryless closed-loop case, (\ref{cbq}), is the same as that, (\ref{q}), in the open-loop case. Thus, we have $q^{**}_B=q^{*}_B$ given the value of $m_i(t)$. 

Since $\frac{\partial q_i}{\partial p_j(t)}>0$, $\gamma''(k_i(t))>0$, $\frac{\partial C(m_i(t),q_i(t))}{\partial m_i(t)}>0$ and $\frac{\partial \Gamma(k^{**}_B,(n-1)k^{**}_B)}{\partial k_i(t)}>0$, if the prices of goods of the firms are strategic substitutes ($\frac{\partial p_j(t)}{\partial m_i(t)}<0$), we have $k^{**}_B>k^{*}_B$; and if the prices of goods of the firms are strategic complements ($\frac{\partial p_j(t)}{\partial m_i(t)}>0$), we have $k^{**}_B<k^{*}_B$.

Note that in the case of linear demand functions the prices of the goods of the firms are strategic complements, and so the R\&D investment of each firm given the cost level in the memoryless closed-loop case is smaller than that in the open-loop case.

\subsection{Cournot oligopoly}

We seek to the solution of the memoryless closed-loop approach in the Cournot oligopoly. Similarly to the previous case, for simplicity, we assume 
\[\frac{\partial \Gamma_i}{\partial K_{-i}(t)}=0.\]
The first order conditions for Firm $i$ are
\begin{align}
\frac{\partial \hat{\mathcal{H}}_i}{\partial q_i(t)}=&p_i(q_1(t),q_2(t),\dots,q_n(t))+\frac{\partial p_i(q_1(t),q_2(t),\dots,q_n(t))}{\partial q_i(t)}q_i(t)\label{c-cbq}\\
&-\frac{\partial C(m_i(t),q_i(t))}{\partial q_i(t)}=0,\notag
\end{align}
and
\begin{equation}
\frac{\partial \hat{\mathcal{H}}_i}{\partial k_i(t)}=-\gamma'(k_i(t))-\lambda_{ii}(t)\frac{\partial \Gamma_i}{\partial k_i(t)}m_i(t)=0.\label{c-cbk}
\end{equation}
The adjoint conditions are
\begin{align}
&-\frac{\partial \hat{\mathcal{H}}_i}{\partial m_i(t)}-\sum_{j\neq i}\frac{\partial \hat{\mathcal{H}}_i}{\partial k_j(t)}\frac{\partial k_j(t)}{\partial m_i(t)}-\sum_{j\neq i}\frac{\partial \hat{\mathcal{H}}_i}{\partial q_j(t)}\frac{\partial q_j(t)}{\partial m_i(t)}\label{c-a11}\\
=&\frac{\partial \lambda_{ii}(t)}{\partial t}-\rho \lambda_{ii}(t),\ i\in \{1, 2, \dots,n\}, \notag
\end{align}
and
\begin{align*}
&-\frac{\partial \hat{\mathcal{H}}_i}{\partial m_j(t)}-\sum_{l\neq i,j}\frac{\partial \hat{\mathcal{H}}_i}{\partial k_l(t)}\frac{\partial k_l(t)}{\partial m_j(t)}-\frac{\partial \hat{\mathcal{H}}_i}{\partial k_j(t)}\frac{\partial k_j(t)}{\partial m_j(t)}
&-\sum_{l\neq i,j}\frac{\partial \hat{\mathcal{H}}_i}{\partial q_l(t)}\frac{\partial q_l(t)}{\partial m_j(t)}-\frac{\partial \hat{\mathcal{H}}_i}{\partial q_j(t)}\frac{\partial q_j(t)}{\partial m_j(t)}\\
=&\frac{\partial \lambda_{ij}(t)}{\partial t}-\rho \lambda_{ij}(t),\ j\neq i.\notag
\end{align*}
We have
\begin{equation}
\frac{\partial \hat{\mathcal{H}}_i}{\partial m_i(t)}=-\frac{\partial C_i}{\partial m_i(t)}+\lambda_{ii}(t)[-\Gamma(k_i(t),K_{-i}(t))+\delta],\label{c-b11}
\end{equation}
\begin{equation}
\frac{\partial \hat{\mathcal{H}}_i}{\partial m_j(t)}=\lambda_{ij}(t)[-\Gamma(k_j(t),K_{-j}(t))+\delta],\label{c-b21}
\end{equation}
\begin{equation*}
\frac{\partial \hat{\mathcal{H}}_i}{\partial k_j(t)}=-\lambda_{ij}(t)\frac{\partial \Gamma_j}{\partial k_j(t)}m_j(t),
\end{equation*}
\begin{equation*}
\frac{\partial k_j(t)}{\partial m_i(t)}=0,
\end{equation*}
\begin{equation*}
\frac{\partial k_j(t)}{\partial m_j(t)}=-\frac{\lambda_{jj}(t)}{\gamma''(k_j(t))+\lambda_{jj}(t)\frac{\partial^2 \Gamma_j}{\partial k_j(t)^2}m_i(t)}\frac{\partial \Gamma_j}{\partial k_j(t)},
\end{equation*}
and
\[\frac{\partial \hat{\mathcal{H}}_i}{\partial q_j(t)}=\frac{\partial p_i}{\partial q_j(t)}q_i(t).\]
$\frac{\partial q_j(t)}{\partial m_i(t)}$ is obtained by (\ref{2-ap3}) in Appendix 3. If the outputs of the firms are strategic substitutes $\left(\frac{\partial p_j}{\partial q_i(t)}+\frac{\partial^2 p_j}{\partial q_i(t)\partial q_j(t)}q_j(t)<0\right)$, $\frac{\partial q_j(t)}{\partial m_i(t)}>0$, and if they are strategic complements $\left(\frac{\partial p_j}{\partial q_i(t)}+\frac{\partial^2 p_j}{\partial q_i(t)\partial q_j(t)}q_j(t)<0\right)$, $\frac{\partial q_j(t)}{\partial m_i(t)}<0$.

At the steady state we have
\[\frac{dm_i(t)}{dt}=m_i(t)[-\Gamma(k_i(t), K_{-i}(t))+\delta]=0,\]
$\frac{\partial \lambda_{ii}}{\partial t}=0$ and $\frac{\partial \lambda_{ij}}{\partial t}=0,\ i\in \{1, 2, \dots, n\},\ j\neq i$. Then, (\ref{c-b11}) and (\ref{c-b21}) are reduced to
\begin{equation*}
\frac{\partial \hat{\mathcal{H}}_i}{\partial m_i(t)}=-\frac{\partial C_i}{\partial m_i(t)},
\end{equation*}
\begin{equation*}
\frac{\partial \hat{\mathcal{H}}_i}{\partial m_j(t)}=0.
\end{equation*}
By symmetry of the oligopoly we can assume $\lambda_{ii}(t)=\lambda_{jj}(t)$ for $j\neq i$, $\lambda_{ij}(t)=\lambda_{il}(t)=\lambda_{ji}(t)$ for $j, l\neq i$, $\frac{\partial k_j(t)}{\partial m_i(t)}=\frac{\partial k_l(t)}{\partial m_j(t)}$ for $l\neq j$, $m_i(t)=m_j(t)$ for $j\neq i$, $q_j(t)=q_i(t)$, and so on.  Denote the steady state values of $\lambda_{ii}$, $\lambda_{ij}$, $q_i(t)$ and $k_i(t)$ by $\lambda_{own}$, $\lambda_{other}$, $q^{**}_C$ and $k^{**}_C$. Then, (\ref{c-a11}) is rewritten as
\[\frac{\partial C_i}{\partial m_i(t)}-(n-1)\frac{\partial p_i}{\partial q_j(t)}q^{**}_C\frac{\partial q_j(t)}{\partial m_i(t)}=-\rho \lambda_{own}.\]
The first order condition for the choice of $k_i(t)$, (\ref{c-cbk}), is reduced to
\[-\gamma'(k^{**}_C)-\lambda_{own}\frac{\partial \Gamma_i}{\partial k_i(t)}m_i(t)=0.\]
This means
\begin{align}
&\left[\frac{\partial C(m_i(t),q_i(t))}{\partial m_i(t)}-(n-1)\frac{\partial p_i}{\partial q_j(t)}q^{**}_C\frac{\partial q_j(t)}{\partial m_i(t)}\right]\frac{\partial \Gamma(k^{**}_C,(n-1)k^{**}_C)}{\partial k_i(t)}\label{c-x2}\\
=&\frac{\rho}{m_i(t)}\gamma'(k^{**}_C).\notag
\end{align}
The first order condition for the output choice in the memoryless closed-loop case, (\ref{c-cbq}), is the same as that, (\ref{cq}), in the open-loop case. Thus, we have $q^{**}_C=q^{*}_C$ given the value of $m_i(t)$. 

Since $\frac{\partial p_i}{\partial q_j(t)}<0$, $\gamma''(k_i(t))>0$, $\frac{\partial C(m_i(t),q_i(t))}{\partial m_i(t)}>0$ and $\frac{\partial \Gamma(k^{**}_C,(n-1)k^{**}_C)}{\partial k_i(t)}>0$, if the outputs of the firms are strategic substitutes ($\frac{\partial q_j(t)}{\partial m_i(t)}>0$), we have $k^{**}_C>k^{*}_C$; and if the outputs of the firms are strategic complements ($\frac{\partial q_j(t)}{\partial m_i(t)}<0$), we have $k^{**}_C<k^{*}_C$.

Note that in the case of linear inverse demand functions the outputs of the firms are strategic substitutes, and so the R\&D investment of each firm given the cost level in the memoryless closed-loop case is larger than that in the open-loop case.

Summarizing the results in this section.

\begin{proposition}
\begin{enumerate}
\item Suppose that there is no spillover effect of R\&D investment in the Bertrand oligopoly. If the prices of goods of the firms are strategic substitutes (or strategic complements), the R\&D investment of each firm given the cost level in the memoryless closed-loop case is larger (or smaller) than that in the open-loop case.

\item Suppose that there is no spillover effect of R\&D investment in the Cournot oligopoly. If the outputs of the firms are strategic substitutes (or strategic complements), the R\&D investment of each firm given the cost level in the memoryless closed-loop case is larger (smaller) than that in the open-loop case.

\item If the demand functions are linear in the Bertrand oligopoly, the R\&D investment of each firm given the cost level in the memoryless closed-loop case is smaller than that in the open-loop case.

\item If the inverse demand functions are linear in the Cournot oligopoly, the R\&D investment of each firm given the cost level in the memoryless closed-loop case is larger than that in the open-loop case.

\end{enumerate}
\end{proposition}

\section{Feedback solution  without spillover}

\subsection{Bertrand oligopoly}

We consider a solution of feedback approach in the Bertrand oligopoly using the Hamilton-Jacobi-Bellman (HJB) equation. Similarly to the memoryless closed-loop case, we assume 
\[\frac{\partial \Gamma_i}{\partial K_{-i}(t)}=0,\]
that is, there is no spillover effect of the R\&D investments. Let $V_i(m_1(t),m_2(t),\dots,m_n(t))$ be the value function of Firm $i,\ i\in \{1, 2, \dots, n\}$. The HJB equation for Firm $i$ is written as
\begin{align}
\rho &V_i(m_1(t),m_2(t),\dots,m_n(t))\label{hjbb}\\
&=\max_{p_i(t),k_i(t)}\{q(p_1(t),p_2(t),\dots,p_n(t))p_i(t)-C(m_i(t),q_i(t))-\gamma_i(k_i(t))\notag\\
&+\frac{\partial V_i(m_1(t),m_2(t),\dots,m_n(t)}{\partial m_i(t)}m_i(t)(-\Gamma_i+\delta)\notag\\
&+\sum_{j\neq i}\frac{\partial V_i(m_1(t),m_2(t),\dots,m_n(t)}{\partial m_j(t)}m_j(t)(-\Gamma_j+\delta)\}.\notag
\end{align}
The first order conditions are
\begin{align}
\frac{\partial \hat{\mathcal{H}}_i}{\partial p_i(t)}=&q_i(p_1(t),p_2(t),\dots,p_n(t)) \label{f1b}\\
&+\left(p_i(t)-\frac{\partial C(m_i(t),q_i(t))}{\partial q_i(t)}\right)\frac{\partial q_i(p_1(t),p_2(t),\dots,p_n(t)}{\partial p_i(t)}=0,\notag
\end{align}
and
\begin{equation}
-\gamma'(k_i(t))-\frac{\partial V_i(m_1(t),m_2(t),\dots,m_n(t))}{\partial m_i(t)}m_i(t)\frac{\partial \Gamma_i}{\partial k_i(t)}=0.\label{f2b}
\end{equation}
From this
\begin{equation*}
\frac{\partial V_i(m_1(t),m_2(t),\dots,m_n(t))}{\partial m_i(t)}m_i(t)=-\frac{\gamma'(k_i(t))}{\frac{\partial \Gamma_i}{\partial k_i(t)}}.
\end{equation*}
Substituting this into (\ref{hjbb}), using symmetry, yields
\begin{align*}
\rho &V_i(m_1(t),m_2(t),\dots,m_n(t))=q_i(p_1(t),p_2(t),\dots,p_n(t))p_i(t)-C(m_i(t),q_i(t))-\gamma_i(k_i(t))\\
&-\frac{\gamma'(k_i(t))}{\frac{\partial \Gamma_i}{\partial k_i(t)}}(-\Gamma_i+\delta)+(n-1)\frac{\partial V_i(m_1(t),m_2(t),\dots,m_n(t)}{\partial m_j(t)}m_j(t)(-\Gamma_j+\delta),\ j\neq i.\notag
\end{align*}
This is an identity. Differentiating it with respect to $m_i(t)$ yields
\begin{align*}
\rho &\frac{\partial V_i(m_1(t),m_2(t),\dots,m_n(t))}{\partial m_i(t)}=-\frac{\partial C(m_i(t),q_i(t))}{\partial m_i(t)}\\
&+\left[q_i(p_1(t),p_2(t),\dots,p_n(t))+\left(p_i(t)-\frac{\partial C(m_i(t),q_i(t))}{\partial q_i(t)}\right)\frac{\partial q_i(p_1(t),p_2(t),\dots,p_n(t)}{\partial p_i(t)}\right]\frac{\partial p_i(t)}{\partial m_i(t)}\\
&+(n-1)\frac{\partial q_j}{\partial p_i(t)}p_i(t)\frac{\partial p_j(t)}{\partial m_i(t)}-\frac{\partial}{\partial k_i(t)}\left(\frac{\gamma'(k_i(t))}{\frac{\partial \Gamma_i}{\partial k_i(t)}}\right)\frac{\partial k_i(t)}{\partial m_i(t)}(-\Gamma_i+\delta)\\
&+(n-1)\frac{\partial^2 V_i(m_1(t),m_2(t),\dots,m_n(t)}{\partial m_i(t)\partial m_j(t)}m_j(t)(-\Gamma_j+\delta).
\end{align*}
At the steady state $-\Gamma_i+\delta=-\Gamma_j+\delta=0$. Thus, using (\ref{f1b}), we get
\begin{align}
\rho \frac{\partial V_i(m_1(t),m_2(t),\dots,m_n(t))}{\partial m_i(t)}=&-\frac{\partial C(m_i(t),q_i(t))}{\partial m_i(t)}\label{f3b}\\
&+(n-1)\frac{\partial q_j}{\partial p_i(t)}p_i(t)\frac{\partial p_j(t)}{\partial m_i(t)}.\notag
\end{align}
From (\ref{f2b}) and (\ref{f3b}), we obtain
\[-\frac{\rho}{m_i(t)}\gamma'(k_i(t))+\frac{\partial \Gamma_i}{\partial k_i(t)}\left[\frac{\partial C(m_i(t),q_i(t))}{\partial m_i(t)}-(n-1)\frac{\partial q_j}{\partial p_i(t)}p_i(t)\frac{\partial p_j(t)}{\partial m_i(t)}\right]=0.\]
This is the same as (\ref{x2}) in the memoryless closed-loop case. Therefore, we get the following proposition.
\begin{proposition}
If there is no  spillover effect of R\&D investment, the memoryless closed-loop solution and the feedback solution  in the Bertrand oligopoly are equivalent.
\end{proposition}

\subsection{Cournot oligopoly}

We consider a solution of feedback approach in the Cournot oligopoly using the HJB equation. Similarly to the previous section we assume 
\[\frac{\partial \Gamma_i}{\partial K_{-i}(t)}=0.\]
Let $V_i(m_1(t),m_2(t),\dots,m_n(t))$ be the value function of Firm $i,\ i\in \{1, 2, \dots, n\}$. The HJB equation for Firm $i$ is written as
\begin{align}
&\rho V_i(m_1(t),m_2(t),\dots,m_n(t))\label{hjb}\\
&=\max_{q_i(t),k_i(t)}\{p_i(q_1(t),q_2(t),\dots,q_n(t))q_i(t)-C(m_i(t),q_i(t))-\gamma_i(k_i(t))\notag\\
&+\frac{\partial V_i(m_1(t),m_2(t),\dots,m_n(t)}{\partial m_i(t)}m_i(t)(-\Gamma_i+\delta)\notag \\
&+\sum_{j\neq i}\frac{\partial V_i(m_1(t),m_2(t),\dots,m_n(t)}{\partial m_j(t)}m_j(t)(-\Gamma_j+\delta)\}.\notag
\end{align}
The first order conditions are
\begin{align}
&p_i(q_1(t),q_2(t),\dots,q_n(t))+\frac{\partial p_i(q_1(t),q_2(t),\dots,q_n(t))}{\partial q_i(t)}q_i(t)\label{f1}\\
&-\frac{\partial C(m_i(t),q_i(t))}{\partial q_i(t)}=0,\notag
\end{align}
and
\begin{equation}
-\gamma'(k_i(t))-\frac{\partial V_i(m_1(t),m_2(t),\dots,m_n(t))}{\partial m_i(t)}m_i(t)\frac{\partial \Gamma_i}{\partial k_i(t)}=0.\label{f2}
\end{equation}
From this
\begin{equation*}
\frac{\partial V_i(m_1(t),m_2(t),\dots,m_n(t))}{\partial m_i(t)}m_i(t)=-\frac{\gamma'(k_i(t))}{\frac{\partial \Gamma_i}{\partial k_i(t)}}.
\end{equation*}
Substituting this into (\ref{hjb}), using symmetry, yields
\begin{align*}
\rho &V_i(m_1(t),m_2(t),\dots,m_n(t))=p_i(q_1(t),q_2(t),\dots,q_n(t))q_i(t)-C(m_i(t),q_i(t))-\gamma_i(k_i(t))\\
&-\frac{\gamma'(k_i(t))}{\frac{\partial \Gamma_i}{\partial k_i(t)}}(-\Gamma_i+\delta)+(n-1)\frac{\partial V_i(m_1(t),m_2(t),\dots,m_n(t)}{\partial m_j(t)}m_j(t)(-\Gamma_j+\delta),\ j\neq i.\notag
\end{align*}
This is an identity. Differentiating it with respect to $m_i(t)$ yields
\begin{align*}
\rho &\frac{\partial V_i(m_1(t),m_2(t),\dots,m_n(t))}{\partial m_i(t)}=-\frac{\partial C(m_i(t),q_i(t))}{\partial m_i(t)}\\
&+\left[p_i(q_1(t),q_2(t),\dots,q_n(t))+\frac{\partial p_i(q_1(t),q_2(t),\dots,q_n(t))}{\partial q_i(t)}q_i(t)-\frac{\partial C(m_i(t),q_i(t))}{\partial q_i(t)}\right]\frac{\partial q_i(t)}{\partial m_i(t)}\\
&+(n-1)\frac{\partial p_i}{\partial q_j(t)}q_i(t)\frac{\partial q_j(t)}{\partial m_i(t)}-\frac{\partial}{\partial k_i(t)}\left(\frac{\gamma'(k_i(t))}{\frac{\partial \Gamma_i}{\partial k_i(t)}}\right)\frac{\partial k_i(t)}{\partial m_i(t)}(-\Gamma_i+\delta)\\
&+(n-1)\frac{\partial^2 V_i(m_1(t),m_2(t),\dots,m_n(t)}{\partial m_i(t)\partial m_j(t)}m_j(t)(-\Gamma_j+\delta).
\end{align*}
At the steady state $-\Gamma_i+\delta=-\Gamma_j+\delta=0$. Thus, using (\ref{f1}), we get
\begin{align}
\rho \frac{\partial V_i(m_1(t),m_2(t),\dots,m_n(t))}{\partial m_i(t)}=&-\frac{\partial C(m_i(t),q_i(t))}{\partial m_i(t)}\label{f3}\\
&+(n-1)\frac{\partial p_i}{\partial q_j(t)}q_i(t)\frac{\partial q_j(t)}{\partial m_i(t)}.\notag
\end{align}
From (\ref{f2}) and (\ref{f3}), we obtain
\[-\frac{\rho}{m_i(t)}\gamma'(k_i(t))+\frac{\partial \Gamma_i}{\partial k_i(t)}\left[\frac{\partial C(m_i(t),q_i(t))}{\partial m_i(t)}-(n-1)\frac{\partial p_i}{\partial q_j(t)}q_i(t)\frac{\partial q_j(t)}{\partial m_i(t)}\right]=0.\]
This is the same as (\ref{c-x2}) in the memoryless closed-loop case. Therefore, we get the following proposition.
\begin{proposition}
If there is no  spillover effect of R\&D investment, the memoryless closed-loop solution and the feedback solution  in the Cournot oligopoly are equivalent.
\end{proposition}

\section{Concluding Remark}

In this paper we analyzed the memoryless closed-loop solution and the feedback solution only when there is no spillover effect of R\&D investment. In the future research we want to investigate the relations among the open-loop solution, the closed-loop solution and the feedback solution in a case with spillovers.


\section*{Acknowledgment}

This work was supported by Japan Society for the Promotion of Science KAKENHI Grant Number 18K01594.

\bibliographystyle{plainnat} 
\bibliography{yatanaka}

\begin{thebibliography}{10}
\providecommand{\natexlab}[1]{#1}
\providecommand{\url}[1]{\texttt{#1}}
\expandafter\ifx\csname urlstyle\endcsname\relax
  \providecommand{\doi}[1]{doi: #1}\else
  \providecommand{\doi}{doi: \begingroup \urlstyle{rm}\Url}\fi

\bibitem[Cellini and Lambertini(2003)]{cl0}
R.~Cellini and L.~Lambertini.
\newblock Advertising in a differential oligopoly game.
\newblock \emph{Journal of Optimization Theory and Applications}, 116:\penalty0
  61--81, 2003.

\bibitem[Cellini and Lambertini(2004)]{cl1}
R.~Cellini and L.~Lambertini.
\newblock Dynamic oligopoly with sticky prices: Closed-loop, feedback and
  open-loop solutions.
\newblock \emph{Journal of Dynamical and Control Systems}, 10:\penalty0
  303--314, 2004.

\bibitem[Cellini and Lambertini(2005)]{cl3}
R.~Cellini and L.~Lambertini.
\newblock {R}\&{D} incentives and market structure: Dynamic analysis.
\newblock \emph{Journal of Optimization Theory and Applications}, 126:\penalty0
  85--96, 2005.

\bibitem[Cellini and Lambertini(2007)]{cl2}
R.~Cellini and L.~Lambertini.
\newblock A differential oligopoly game with differentiated goods and sticky
  prices.
\newblock \emph{European Journal of Operational Research}, 176:\penalty0
  1131--1144, 2007.

\bibitem[Cellini and Lambertini(2009)]{cl5}
R.~Cellini and L.~Lambertini.
\newblock Dynamic {R}\&{D} with spillovers: Competition and control.
\newblock \emph{Journal of Economic Dynamics and Control}, 33:\penalty0
  568--582, 2009.

\bibitem[Cellini and Lambertini(2011)]{cl4}
R.~Cellini and L.~Lambertini.
\newblock {R}\&{D} incentives under {B}ertrand competition: A differential
  game.
\newblock \emph{Japanese Economic Review}, 62:\penalty0 387--400, 2011.

\bibitem[Fujiwara(2006)]{fu}
Kenji Fujiwara.
\newblock A stackelberg game model of dynamic duopolistic competition with
  sticky prices.
\newblock \emph{Economics Bulletin}, 12:\penalty0 1--9, 2006.

\bibitem[Fujiwara(2008)]{fu1}
Kenji Fujiwara.
\newblock Duopoly can be more anti-competitive than monopoly.
\newblock \emph{Economics Letters}, 101:\penalty0 217--219, 2008.

\bibitem[Lambertini(2018)]{lam18}
L.~Lambertini.
\newblock \emph{Differential Games in Industrial Economics}.
\newblock Cambridge University Press, 2018.

\bibitem[Smrkolj and Wagener(2016)]{smr}
G.~Smrkolj and F.~Wagener.
\newblock Dynamic {R}\&{D} with splillovers: A comment.
\newblock \emph{Journal of Economic Dynamics and Control}, 73:\penalty0
  453--457, 2016.

\end{thebibliography}

\appendix

\normalsize

\section*{\boldmath  Appendix 1: Derivation of (\ref{ee1})}

The direct and the inverse demand functions are as follows.
\begin{equation*}
q_i=q_i(p_1, p_2, \dots, p_n),\ i\in \{1, 2, \dots, n\},
\end{equation*}
\begin{equation}
p_i=p_i(q_1, q_2, \dots, q_n),\ i\in \{1, 2, \dots, n\}.\label{e0}
\end{equation}
We omit $t$. Differentiating (\ref{e0}) with respect to $p_i$ given $p_j,\ j\in \{1,2,\dots,n\},\ j\neq i$, yields
\[\frac{\partial p_i}{\partial q_i}\frac{\partial q_i}{\partial p_i}+\sum_{j=1, j\neq i}^n\frac{\partial p_i}{\partial q_j}\frac{\partial q_j}{\partial p_i}=1,\]
and
\[\frac{\partial p_j}{\partial q_i}\frac{\partial q_i}{\partial p_i}+\frac{\partial p_j}{\partial q_j}\frac{\partial q_j}{\partial p_i}+\sum_{k=1, k\neq i,j}^n\frac{\partial p_j}{\partial q_k}\frac{\partial q_k}{\partial p_i}=0,\ j\in \{1, 2, \dots, n\},\ j\neq i.\]
Since by symmetry $\frac{\partial p_i}{\partial q_j}=\frac{\partial p_j}{\partial q_i}=\frac{\partial p_j}{\partial q_k}$, $\frac{\partial p_j}{\partial q_j}=\frac{\partial p_i}{\partial q_i}$ and $\frac{\partial q_k}{\partial p_i}=\frac{\partial q_j}{\partial p_i}$ at the steady state, they are rewritten as
\[\frac{\partial p_i}{\partial q_i}\frac{\partial q_i}{\partial p_i}+(n-1)\frac{\partial p_j}{\partial q_i}\frac{\partial q_j}{\partial p_i}=1,\]
and
\[\frac{\partial p_j}{\partial q_i}\frac{\partial q_i}{\partial p_i}+\left[\frac{\partial p_i}{\partial q_i}+(n-2)\frac{\partial p_j}{\partial q_i}\right]\frac{\partial q_j}{\partial p_i}=0.\]
From them we get
\begin{equation}
\frac{\partial q_i}{\partial p_i}=\frac{\partial q_j}{\partial p_j}=\frac{\frac{\partial p_i}{\partial q_i}+(n-2)\frac{\partial p_j}{\partial q_i}}{\left(\frac{\partial p_i}{\partial q_i}-\frac{\partial p_j}{\partial q_i}\right)\left[\frac{\partial p_i}{\partial q_i}+(n-1)\frac{\partial p_j}{\partial q_i}\right]},\label{e11}
\end{equation}
and
\begin{equation*}
\frac{\partial q_j}{\partial p_i}=\frac{\partial q_i}{\partial p_j}=-\frac{\frac{\partial p_j}{\partial q_i}}{\left(\frac{\partial p_i}{\partial q_i}-\frac{\partial p_j}{\partial q_i}\right)\left[\frac{\partial p_i}{\partial q_i}+(n-1)\frac{\partial p_j}{\partial q_i}\right]},
\end{equation*}
because $\frac{\partial q_i}{\partial p_j}=\frac{\partial q_j}{\partial p_i}$ and $\frac{\partial q_i}{\partial p_i}=\frac{\partial q_j}{\partial p_j}$ at the steady state.

From (\ref{e11})
\begin{align}
&1-\frac{\partial q_i}{\partial p_i}\frac{\partial p_i}{\partial q_i}=\frac{\left(\frac{\partial p_i}{\partial q_i}-\frac{\partial p_j}{\partial q_i}\right)\left[\frac{\partial p_i}{\partial q_i}+(n-1)\frac{\partial p_j}{\partial q_i}\right]-\left(\frac{\partial p_i}{\partial q_i}\right)^2-(n-2)\frac{\partial p_i}{\partial q_i}\frac{\partial p_j}{\partial q_i}}{\left(\frac{\partial p_i}{\partial q_i}-\frac{\partial p_j}{\partial q_i}\right)\left[\frac{\partial p_i}{\partial q_i}+(n-1)\frac{\partial p_j}{\partial q_i}\right]}\label{1-a}\\
=&\frac{-(n-1)\left(\frac{\partial p_j}{\partial q_i}\right)^2}{\left(\frac{\partial p_i}{\partial q_i}-\frac{\partial p_j}{\partial q_i}\right)\left[\frac{\partial p_i}{\partial q_i}+(n-1)\frac{\partial p_j}{\partial q_i}\right]}<0.\notag
\end{align}

\section*{\boldmath Appendix 2: Derivation of $\frac{\partial p_j(t)}{\partial m_i(t)}$.}

Suppose a state such that $p_1(t)=p_2(t)=\dots=p_n(t)$. The first order conditions for Firm $i$ and Firm $j,\ j\neq i$, are
\begin{align*}
q_i(p_1(t),p_2(t),\dots,p_n(t))+\left(p_i(t)-\frac{\partial C_i}{\partial q_i(t)}\right)\frac{\partial q_i(p_1(t),p_2(t),\dots,p_n(t)}{\partial p_i(t)}=0,\notag
\end{align*}
and
\begin{align*}
q_j(p_1(t),p_2(t),\dots,p_n(t))+\left(p_j(t)-\frac{\partial C_j}{\partial q_j(t)}\right)\frac{\partial q_j(p_1(t),p_2(t),\dots,p_n(t)}{\partial p_j(t)}=0,\notag
\end{align*}
Denote $q_i(p_1(t), p_2(t), \dots, p_n(t))$ by $q_i$. Differentiating them with respect to $m_i(t)$ yields
\[\varphi_i\frac{\partial p_i(t)}{\partial m_i(t)}+\psi_i\frac{\partial p_j(t)}{\partial m_i(t)}=\frac{\partial^2 C_i}{\partial q_i(t)\partial m_i(t)}\frac{\partial q_i}{\partial p_i(t)},\]
and
\[\psi_j\frac{\partial p_i(t)}{\partial m_i(t)}+\varphi_j\frac{\partial p_j(t)}{\partial m_i(t)}=0,\]
where
\[\varphi_i=2\frac{\partial q_i}{\partial p_i(t)}+\left(p_i(t)-\frac{\partial C_i}{\partial q_i(t)}\right)\frac{\partial^2 q_i}{\partial p_i(t)^2}-\frac{\partial^2 C_i}{\partial q_i(t)^2}\left(\frac{\partial q_i}{\partial p_i(t)}\right)^2,\]
\[\psi_i=(n-1)\left[\frac{\partial q_i}{\partial p_j(t)}+\left(p_i(t)-\frac{\partial C_i}{\partial q_i(t)}\right)\frac{\partial^2 q_i}{\partial p_i(t)\partial p_j(t)}\right],\]
\[\psi_j=\frac{\partial q_j}{\partial p_i(t)}+\left(p_j(t)-\frac{\partial C_j}{\partial q_j(t)}\right)\frac{\partial^2 q_j}{\partial p_i(t)\partial p_j(t)},\]
\begin{align*}
\varphi_j=&2\frac{\partial q_j}{\partial p_j(t)}+\left(p_j(t)-\frac{\partial C_j}{\partial q_j(t)}\right)\frac{\partial^2 q_j}{\partial p_j(t)^2}\\
&+(n-2)\left[\frac{\partial q_j}{\partial p_l(t)}+\left(p_j(t)-\frac{\partial C_j}{\partial q_j(t)}\right)\frac{\partial^2 q_j}{\partial p_j(t)\partial p_l(t)}\right]\\
&-\frac{\partial^2 C_j}{\partial q_j(t)^2}\left(\frac{\partial q_j}{\partial p_j(t)}\right)^2,\ l\neq j.
\end{align*}
From them we obtain
\[\frac{\partial p_i(t)}{\partial m_i(t)}=\frac{\varphi_j}{\varphi_i\varphi_j-\psi_i\psi_j}\frac{\partial^2 C_i}{\partial q_i(t)\partial m_i(t)}\frac{\partial q_i(t)}{\partial p_i(t)}<0,\]
and
\begin{equation}
\frac{\partial p_j(t)}{\partial m_i(t)}=-\frac{\psi_j}{\varphi_i\varphi_j-\psi_i\psi_j}\frac{\partial^2 C_i}{\partial q_i(t)\partial m_i(t)}\frac{\partial q_i(t)}{\partial p_i(t)}.\label{2ap3}
\end{equation}
We have $\varphi_i<0$, $\varphi_j<0$ and $\varphi_i\varphi_j-\psi_i\psi_j>0$. If the prices of goods of the firms are strategic substitutes $\left(\psi_j<0\right)$, $\frac{\partial p_j(t)}{\partial m_i(t)}<0$, and if they are strategic complements $\left(\psi_j>0\right)$, $\frac{\partial p_j(t)}{\partial m_i(t)}>0$.

\section*{Appendix 3: Derivation of $\frac{\partial q_j(t)}{\partial m_i(t)}$.}

Suppose a state such that $q_1(t)=q_2(t)=\dots=q_n(t)$. The first order conditions for Firm $i$ and Firm $j,\ j\neq i$ are
\begin{equation*}
p_i(q_1(t), q_1(t), \dots, q_n(t))+\frac{\partial p_i(q_1(t), q_1(t), \dots, q_n(t))}{\partial q_i(t)}q_i(t)-\frac{\partial C_i}{\partial q_i(t)}=0,
\end{equation*}
and
\begin{equation*}
p_j(q_1(t), q_1(t), \dots, q_n(t))+\frac{\partial p_j(q_1(t), q_1(t), \dots, q_n(t))}{\partial q_j(t)}q_j(t)-\frac{\partial C_j}{\partial q_j(t)}=0,
\end{equation*}
Denote $p_i(q_1(t), q_1(t), \dots, q_n(t))$ by $p_i$. Differentiating them with respect to $m_i(t)$ yields
\begin{align*}
&\left[2\frac{\partial p_i}{\partial q_i(t)}+\frac{\partial^2 p_i}{\partial q_i(t)^2}q_i(t)-\frac{\partial^2 C_i}{\partial q_i(t)^2}\right]\frac{\partial q_i(t)}{\partial m_i(t)}+(n-1)\left[\frac{\partial p_i}{\partial q_j(t)}+\frac{\partial^2 p_i}{\partial q_i(t)\partial q_j(t)}q_i(t)\right]\frac{\partial q_j(t)}{\partial m_i(t)}\\
=&\frac{\partial^2 C_i}{\partial q_i(t)\partial m_i(t)},
\end{align*}
and
\[\left[\frac{\partial p_j}{\partial q_i(t)}+\frac{\partial^2 p_j}{\partial q_i(t)\partial q_j(t)}q_j(t)\right]\frac{\partial q_i(t)}{\partial m_i(t)}+\left[n\frac{\partial p_j}{\partial q_j(t)}+(n-1)\frac{\partial^2 p_j}{\partial q_j(t)^2}q_j(t)-\frac{\partial^2 C_j}{\partial q_j(t)^2}\right]\frac{\partial q_j(t)}{\partial m_i(t)}=0.\]
From them we obtain
\[\frac{\partial q_i(t)}{\partial m_i(t)}=\frac{n\frac{\partial p_j}{\partial q_j(t)}+(n-1)\frac{\partial^2 p_j}{\partial q_j(t)^2}q_j(t)-\frac{\partial^2 C_j}{\partial q_j(t)^2}}{\Delta}\frac{\partial^2 C_i}{\partial q_i(t)\partial m_i(t)}<0,\]
and
\begin{equation}
\frac{\partial q_j(t)}{\partial m_i(t)}=-\frac{\frac{\partial p_j}{\partial q_i(t)}+\frac{\partial^2 p_j}{\partial q_i(t)\partial q_j(t)}q_j(t)}{\Delta}\frac{\partial^2 C_i}{\partial q_i(t)\partial m_i(t)},\label{2-ap3}
\end{equation}
where
\begin{align*}
\Delta=&\left[2\frac{\partial p_i}{\partial q_i(t)}+\frac{\partial^2 p_i}{\partial q_i(t)^2}q_i(t)-\frac{\partial^2 C_i}{\partial q_i(t)^2}\right]\left[n\frac{\partial p_j}{\partial q_j(t)}+(n-1)\frac{\partial^2 p_j}{\partial q_j(t)^2}q_j(t)-\frac{\partial^2 C_j}{\partial q_j(t)^2}\right]\\
&-(n-1)\left[\frac{\partial p_i}{\partial q_j(t)}+\frac{\partial^2 p_i}{\partial q_i(t)\partial q_j(t)}q_i(t)\right]\left[\frac{\partial p_j}{\partial q_i(t)}+\frac{\partial^2 p_j}{\partial q_i(t)\partial q_j(t)}q_j(t)\right]>0.
\end{align*}
If the outputs of the firms are strategic substitutes $\left(\frac{\partial p_j}{\partial q_i(t)}+\frac{\partial^2 p_j}{\partial q_i(t)\partial q_j(t)}q_j(t)<0\right)$, $\frac{\partial q_j(t)}{\partial m_i(t)}>0$, and if they are strategic complements $\left(\frac{\partial p_j}{\partial q_i(t)}+\frac{\partial^2 p_j}{\partial q_i(t)\partial q_j(t)}q_j(t)>0\right)$, $\frac{\partial q_j(t)}{\partial m_i(t)}<0$.

\end{document}